\input amstex
\TagsOnRight 
\raggedbottom
\magnification=\magstep1         

\font \sec=cmbx12
\font \cd=cmr9 
\font \cmmi=cmmi12 
\baselineskip=15pt

\def\ref#1{\noindent\item{[#1.]}}

\vskip 9cm
\hfill {~~~}
\vskip 2cm
\sec
\centerline {\sec {\cmmi q}-Deformation of the Krichever-Novikov Algebra }
\rm
\vskip.7cm
\centerline {C.H. Oh $^{\dag}$ and K. Singh }
\vskip .5cm
{\sl {\centerline {Department of Physics}}}
{\sl {\centerline {Faculty of Science}}}
{\sl {\centerline {National University of Singapore}}}
{\sl {\centerline {Lower Kent Ridge, Singapore 0511}}}
{\sl {\centerline {Republic of Singapore}}}
\vskip 2cm
{\bf {\centerline {Abstract}}}
\vskip .5cm
\noindent
Using $q$-operator product expansions between $U(1)$ current fields
and also the corresponding energy-momentum tensors, we furnish the
$q$-analoques of the generalized Heisenberg algebra and the
Krichever-Novikov algebra.  
\vskip.5cm
\vskip 2cm 
{\bf {\centerline {To appear in Letters in Mathematical Physics}}}
\vskip 3.5truecm
\hrule width6cm
\vskip .3truecm
{\sevenrm $^{\dag}$ E-mail: PHYOHCH\@NUS.SG}
\eject
\rm
\par
The recent focus on deformations of algebras called quantum
algebras can be attributed to the fact that they appear to be the
basic algebraic structures underlying an amazingly diverse set of
physical situations. To date many interesting features of these
algebras have been found and they are now known to belong to a
class of algebras called Hopf algebras [1]. The remarkable aspect
of these structures is that they can be regarded as deformations of
the usual Lie algebras. Of late, there has been a considerable
interest in the deformation of the Virasoro algebra and the
underlying Heisenberg algebra [2-11]. In this letter we focus our
attention on deforming generalizations of these algebras, namely
the Krichever-Novikov (KN) algebra and its associated Heisenberg
algebra.
\par
The introduction by Krichever and Novikov [12-14] of a basis for the
general tensor field over Riemann surfaces of arbitrary genus 
has led to a new operator formalism which generalizes conformal
field theory on these surfaces. Indeed, on a suface $\Sigma_g$ of
genus $g$, the relevant algebra to consider is the centrally
extended KN algebra:
$$
[{\Cal L}_m,{\Cal L}_n] = \sum_{s=m+n-g_0}^{m+n+g_0} c_{m,n}^s {\Cal L}_{m+n-s} +
c\chi_{mn}\tag 1
$$ 
where $g_0=3g/2$. Here the indices $m,n$ take integral values if $g$ is even 
and half integral values for odd $g$. 
The structure constants $c_{m,n}^{s}$ and the cocycle $\chi_{mn}$ are evaluated
using the KN bases for the vector fields $\{e_n\}$ and the the
quadratic differentials $\{\Omega_n\}$. Specifically one can
introduce local coordinates $z_+$ and $z_-$ around two
distinguished points $P_+$ and $P_-$ with the basis elements 
given by
$$\alignat2
e_n &\equiv e_n(z_{\pm}){{\partial}\over{\partial z_{\pm}}} &&= \sum_{m=0}^{\infty}e_{n,m}^{\pm}z_{\pm}^{\pm n-g_0+1+m}\frac {\partial}
{\partial {z_{\pm}}} \tag 2a\\
\Omega_n &\equiv \Omega_n(z_{\pm})dz_{\pm}^2 &&= \sum_{m=0}^{\infty}\Omega_{n,m}^{\pm}
z_{\pm}^{\mp n+g_0-2+m}(dz_{\pm})^2 \tag 2b
\endalignat$$
in the neighbourhood of $P_{\pm}$. These two bases are dual in the sense
$$\pm \oint_{C_{\pm}}\frac {dz_{\pm}}{2\pi i}e_{m}(z_{\pm})\Omega_{n}(z_{\pm})
= \delta_{mn} \tag 3 $$
where $C_{+} (C_{-})$ denotes any contour around $P_{+} (P_{-})$ but
not including $P_{-} (P_{+})$. In these bases, one has
 $$c_{m,n}^s=\oint_{C_{\pm}}\frac{dw_{\pm}}{2\pi i}
\left({\partial_{w_{\pm}}e_m(w_{\pm})e_n(w_{\pm}) - 
e_m(w_{\pm})\partial_{w_{\pm}}e_n(w_{\pm})}\right)\Omega_{m+n-s}(w_{\pm})
\tag 4a$$
and
$$\chi_{mn}=\frac {1}{12}\oint_{C_{\pm}}\frac{dw_{\pm}}{2\pi i}
e_m'''(w_{\pm})e_n(w_{\pm}).\tag 4b$$
\par
Similarly, the generalized Heisenberg algebra appropriate for $\Sigma_g$
is defined by the relations [13]
$$[\alpha_m,\alpha_n]=\gamma_{mn}. \tag 5$$
The structure constants are evaluated through 
$$\gamma_{mn}=\frac {1}{2\pi i} \oint_{C_{\tau}}dA_mA_n\tag 6$$
where $\{A_n\}$ is a basis of meromorphic functions on $\Sigma_g$.
\par
Now, one of the most remarkable feature of the KN formalism 
is its simplicity in quantizing theories over Riemann surfaces. In fact it 
extends the standard operator formalism
in a most natural way, removing the need for special considerations
arising from the non-trivial topology. Indeed,
the quantization of fields over $\Sigma_g$ is carried out in 
the traditional way in which the coefficients of the fields in a series 
expansion are regarded as operators on some Fock space. Here `globalization' 
of the theory is obtained by replacing the Laurent basis, in which fields are
usually expanded, to these KN bases. The latter intrinsically
carries all the topological information of the underlying surface. In the
field-theoretic context, the full content of the theory is usually
embodied in the operator product expansions {\cd (OPEs)} of the relevant field
operators. For instance, in the case of a conformal field theory,
the {\cd OPE} between two energy-momentum tensors bears all the
information about the underlying algebra, which in the case of a
sphere, is the Virasoro algebra. Although the {\cd OPEs} between field
operators are generically genus dependent, the singular terms 
are not [15]. What is interesting here is that the essential information on the 
algebra is contained only in the singular portion of the {\cd OPE}.  
Using this fact,
Mezincescu et al [16] showed that the {\cd KN} algebra can be obtained from
the {\cd OPE} of $T(z)T(w)$ with the usual Laurent basis replaced
by the {\cd KN} basis.
\plainfootnote {$^{\dag}$}
{\cd Although they consider the torus ($g=1$) exclusively, their
analysis can be extended to a surface of any genus (see ref.[17]).} 
In the following we employ similar techniques to obtain 
$q$-analogues of the generalized Heisenberg algebra and the {\cd KN} algebra. 
\par
To this end, we consider a $q$-deformed version of a $U(1)$ current
algebra and its associated conformal algebra defined on the plane.
By using a $q$-analoque of the Sugawara construction,
$$T^{\alpha}(z)={{1}\over{4}} :J(zq^{\alpha/2})J(zq^{-\alpha/2}): +
{{1}\over{4}} :J(zq^{-\alpha/2})J(zq^{\alpha/2}): \tag 7$$
it was shown in ref.[10] that the $q$-{\cd OPE's}
\plainfootnote {$^{\ddag}$}
{\cd In the expressions, the $q$-product $(z-w)_{q}^{2}$ denotes 
$(z-wq^{-1})(z-wq)$.}  

$$ \alignat 2
J(z)J(w) &\sim {{\kappa}\over
{(z-w)_{q}^{2}}} && {~}  \tag 8\\
T^{\alpha}(z)T^{\beta}(w) 
&\sim {{\kappa}\over {2(q-q^{-1})w}}
\{ &&{{T^{\alpha +\beta +1}(wq^{(\alpha + 1)/2})}\over 
{(zq^{-(\alpha -\beta)/2}-wq^{\beta +1})}} 
+{{T^{-\alpha +\beta -1}(wq^{(\alpha + 1)/2})}\over 
{(zq^{-(\alpha +\beta)/2}-wq^{-\beta +1})}}\\ 
&{~}
&-&{{T^{-\alpha -\beta +1}(wq^{(\alpha - 1)/2})}\over 
{(zq^{-(\alpha -\beta)/2}-wq^{\beta -1})}} 
-{{T^{\alpha -\beta -1}(wq^{(\alpha - 1)/2})}\over 
{(zq^{-(\alpha +\beta)/2}-wq^{-\beta -1})}}\rbrace\\
&+{{\kappa}^2\over {4(q-q^{-1})w^3}}\{
&&{{1}\over{(q^{\alpha +\beta /2 +1}-q^{-\beta /2})_{q}^{2}}}
{{1}\over{(zq^{-(\alpha -\beta) /2 }-wq^{\beta +1})}}\\
&{~}
&+&{{1}\over{(q^{\alpha -\beta /2 +1}-q^{\beta /2})_{q}^{2}}}
{{1}\over{(zq^{-(\alpha +\beta) /2 }-wq^{-\beta +1})}}\\
&{~}
&-&{{1}\over{(q^{\alpha +\beta /2 -1}-q^{-\beta /2})_{q}^{2}}}
{{1}\over{(zq^{-(\alpha -\beta) /2 }-wq^{\beta -1})}}\\
&{~}
&-&{{1}\over{(q^{\alpha -\beta /2 -1}-q^{\beta /2})_{q}^{2}}}
{{1}\over{(zq^{-(\alpha +\beta) /2 }-wq^{-\beta -1})}}\rbrace\\
&{~}
&+&\qquad \qquad q\leftrightarrow q^{-1}\tag 9 \endalignat $$
where $\kappa = (q-q^{-1})/\ln q^2$ 
provide a realization of the $q$-deformed Heisenberg algebra and
the corresponding $q$-deformed Virasoro algebra of Chaichian and
Pre${\check{\text s}}$najder [4].
\par
Now for the higher genus version, we first write the $U(1)$ current
generators $\{\alpha_n\}$ on $\Sigma_g$ canonically as  
$$\alpha_n= {{1}\over{2\pi i}}\oint_{C_{\tau}}J(Q)A_n(Q)\qquad Q\in
\Sigma_g\tag 10$$
where the contour $C_{\tau}$ is an `equal-time' contour [13]. In the
neighbourhood of the distinguished points $P_{\pm}$ with $J(Q)=J(z_{\pm})dz_{\pm}$, 
this takes the form of 
$$\alpha_n= \oint_{C_{\pm}}{{dz_{\pm}}\over{2\pi i}}J(z_{\pm})A_n(z_{\pm})
\tag 11$$
where the contour $C_{\tau}$ is now replaced by $C_{\pm}$. The
commutator betweeen the generators can then be evaluated using the
standard notion of radial ordering. To make this more explicit, we
first consider the undeformed case, we have for $J(z)J(w)\sim (z-w)^{-2}$
$$\align
[\alpha_m,\alpha_n] &= \oint_{C_{\pm}}\frac {dw_{\pm}}{2\pi i} 
\oint_{C_{w}}\frac {dz_{\pm}}{2\pi i} J(z_{\pm})J(w_{\pm}) 
A_m(z_{\pm})A_n(w_{\pm})\\
&= \oint_{C_{\pm}}\frac {dw_{\pm}}{2\pi i} 
\oint_{C_{w}}\frac {dz_{\pm}}{2\pi i} \frac {A_m(z_{\pm})A_n(w_{\pm})}
{(z_{\pm}-w_{\pm})^2}\\
&=\oint_{C_{\pm}}\frac {dw_{\pm}}{2\pi i} {A'}_m(w_{\pm})A_n(w_{\pm})=\gamma_{mn}.
\tag 12 \endalign$$
where $C_w$ is a contour taken around the singular points of the
{\cd OPE} of $J(z)$ and $J(w)$. For the $q$-deformed case,
\plainfootnote {$^{\dag}$}{\cd $q$ is assumed to be real here.}
 we
simply replace the {\cd OPE} by its $q$-analoque (8):   
$$\align
[\alpha_m,\alpha_n] &= \oint_{C_{\pm}}\frac {dw_{\pm}}{2\pi i} 
\oint_{C_{w}}\frac {dz_{\pm}}{2\pi i} \frac {A_m(z_{\pm})A_n(w_{\pm})}
{(z_{\pm}-w_{\pm})_q^2}\\
&=\oint_{C_{\pm}}\frac {dw_{\pm}}{2\pi i} \left({ \partial^{q}_{w_{\pm}}
A_m(w_{\pm})}\right) A_n(w_{\pm})
=\gamma_{mn}^q\tag 13 \endalign $$
where 
$$\partial_{w_{\pm}}^{q}A_m(w_{\pm})\equiv {{A_m(w_{\pm}q)}-A_m(w_{\pm}q^{-1})
\over{w_{\pm}(q-q^{-1})}}.\tag 14$$
It is interesting to note that in the $q\to 1$ limit, one has 
$\partial_{w_{\pm}}^{q}A_m(w_{\pm})\to
\partial_{w_{\pm}}A_m(w_{\pm})$ and consequently $\gamma_{mn}^q\to \gamma_{mn}$.
\par
Next let us consider the $q$-analogue of the Virasoro algebra. By writing
$$\align
T^{\alpha}(z_{\pm})&=\sum_n \Omega_n(z_{\pm}){\Cal L}_{n}^{\alpha}
\tag 15a\\
{\Cal L}^{\alpha}_m &=\pm \oint_{C_{\pm}}\frac {dz_{\pm}}{2\pi
i}e_n(z_{\pm}) T^{\alpha}(z_{\pm})\tag 15b\endalign$$
the commutator between the generators $\{{\Cal L}_m^{\alpha}\}$ can be 
obtained from the $q$-{\cd OPE} (9). To evaluate this, it instructive to
first examine the operator part of the $q$-{\cd OPE} ({\it i.e.} terms
involving the energy-momentum tensor). To this end we consider the
generic term    
$$\frac {\kappa}{2(q-q^{-1})} \frac {T^a(wq^b)}{w(zq^c-wq^d)}\quad
+\quad q\leftrightarrow q^{-1}\tag 16  $$
where the indices $a,b,c,d$ stand for the different algebraic terms
involving $\alpha$ and $\beta$ in the $q$-{\cd OPE} of (9). 
The relevant quantity to be computed here is 
$$\frac {\kappa}{2(q-q^{-1})} \oint_{C_{\pm}}\frac {dw_{\pm}}{2\pi i}
\oint_{C_{w}}\frac {dz_{\pm}}{2\pi i} 
\frac {e_m(z_{\pm})e_n(w_{\pm})T^a(w_{\pm} q^b)}{w(z_{\pm}q^c-w_{\pm}q^d)}
\quad +\quad q\leftrightarrow q^{-1}.\tag 17 $$
By using (15a) and carrying out the $z$-integration this reduces
to 
$$\frac {\kappa}{2(q-q^{-1})}\sum_k  \oint_{C_{\pm}}\frac {dw_{\pm}}{2\pi i}
q^{-c} \frac {e_m(w_{\pm}q^{d-c})e_n(w_{\pm})\Omega_k (w_{\pm} q^b)}{w_{\pm}}
{\Cal L}_k^{a}\quad +\quad q\leftrightarrow q^{-1}.\tag 18 $$
which upon including the $q\leftrightarrow q^{-1}$ term leads to
\plainfootnote{$^{\dag}$}{\cd The $q$-bracket $[x]_q$ denotes ${{q^x-q^{-x}}\over
{q-q^{-1}}}.$}
$$\align
\frac {\kappa}{2}\sum_k \oint_{C_{\pm}}\frac {dw_{\pm}}{2\pi i}
\lbrace &q^{-c}[d-b-c]_q\partial_{w_{\pm}}^{q^{d-b-c}}e_m(w_{\pm})
e_n(w_{\pm}q^{-b})\\
-&q^{-c} [b]_qe_m(w_{\pm}q^{b+c-d})\partial_{w_{\pm}}^{q^{-b}}e_n(w_{\pm})\\
-&[c]_q \frac {e_m(w_{\pm}q^{b+c-d})e_n(w_{\pm} q^{b})}{w_{\pm}}\rbrace 
\Omega_k(w_{\pm}) {\Cal L}^a_k. \tag 19
\endalign$$ 
The sum over $k$ in (19) is not infinite but runs over
$(2g_0+1)$ terms ({\it i.e.} from $m+n-g_0$ to $m+n+g_0$). This can
be easily shown by evaluating the $w$-integral around the points
$P_+$ {\it and} $P_-$ using the explicit form of the basis elements.
Then by denoting
$$\align
{\Cal D}_{m,n}^s(b,c,d)&\equiv
\frac {\kappa}{2}\oint_{C_{\pm}}\frac {dw_{\pm}}{2\pi i}
\lbrace q^{-c}[d-b-c]_q\partial_{w_{\pm}}^{q^{d-b-c}}e_m(w_{\pm})
e_n(w_{\pm}q^{-b})\\
-&q^{-c} [b]_qe_m(w_{\pm}q^{b+c-d})\partial_{w_{\pm}}^{q^{-b}}e_n(w_{\pm})\\
-&[c]_q \frac {e_m(w_{\pm}q^{b+c-d})e_n(w_{\pm} q^{b})}{w_{\pm}}\rbrace 
\Omega_{m+n-s}(w_{\pm}) . \tag 20
\endalign $$   
the operator part of the commutation relations can be written as
$$\alignat2
[{\Cal L}_m,{\Cal L}_n] = \sum_{s=-g_0}^{g_0}
\lbrace &{\Cal D}_{m,n}^s(\frac {\alpha +1}{2}
,\frac {\beta-\alpha }{2},\beta +1) &&{\Cal L}_{m+n-s}^{\alpha +\beta +1}\\
+&{\Cal D}_{m,n}^s(\frac {\alpha +1}{2}
,\frac {-\beta-\alpha }{2},-\beta +1) &&{\Cal L}_{m+n-s}^{-\alpha +\beta +1}\\
-&{\Cal D}_{m,n}^s(\frac {\alpha -1}{2}
,\frac {\beta-\alpha }{2},\beta -1) &&{\Cal L}_{m+n-s}^{-\alpha -\beta +1}\\
-&{\Cal D}_{m,n}^s(\frac {\alpha -1}{2}
,\frac {-\beta-\alpha }{2},-\beta -1) &&{\Cal L}_{m+n-s}^{\alpha -\beta -1}
\rbrace \\
+&\qquad {\text {central term}}. \tag 21 \endalignat$$
Using some of the techniques above, it is not difficult to show
that the central term is given by
$$\chi^{q}_{m,n}=-\frac {1}{16(\ln q)^2}\sum_{k=0}^{2g_0-m-n} e_{m,k}^
+e_{n,2g_0-m-n-k}^+
\left\{{C^{\alpha,\beta}_{m+1}(g_0;k)+C^{\alpha,-\beta}_{m+1}(g_0;k)}\right\}
\tag 22$$
where
$$\align {C}_{m}^{\alpha , \beta}(g_0;k) = &{~}{{[{{\alpha +\beta
+ 2}\over {2}}(m-g_0+k)]_q}\over{{[{{\alpha +\beta +
2}\over {2}} }]_q}}
+ {{[{{\alpha +\beta
- 2}\over {2}} (m-g_0+k)]_q}\over{{[{{\alpha +\beta -
2}\over {2}} }]_q}}\\
&-  (q^{(m-1-g_0+k)}+q^{-(m-1-g_0+k)}) {{[{{\alpha +\beta
}\over {2}} (m-g_0+k)]_q}\over{{[{{\alpha +\beta }\over {2}} 
}]_q}}.\tag 23 \endalign$$
It should be noted that the coefficients $e_{m,k}^+$ can be
evaluated for a given basis of vector fields through
$$e_{m,k}^+=\oint_{C_+}\frac {dz}{2\pi i}e_m(z_+)z^{-m-k+g_0-1}.\tag
24$$  
For instance in the $g=1$ case this can be computed explicitly for
the vector fields $e_n=e_n(z)\partial/\partial z$ with $\{e_n(z)\}$
given in terms of 
of the well studied elliptic functions [18]:
$$\align e_n(z) &=
\frac {\sigma^{n-1/2}(z-z_{0})\sigma(z+2nz_{0})}
{\sigma^{n+1/2}(z+z_{0})}
\frac {\sigma^{n+1/2}(2z_{0})}
{\sigma((2n+1)z_{0})} \qquad n\ne -1/2 \tag 25a \\
e_{-1/2}(z) &=
\frac {\sigma^{2}(z)}
{\sigma(z+z_{0})\sigma(z-z_{0})}
\frac {\sigma(2z_{0})}
{\sigma^2(z_{0})}. \tag 25b \endalign $$
Here $\sigma(z)$ is the Wierstrass sigma-function. It should also
be noted that, by computing $\chi_{m,n}^{q}$ in the vicinity of
$P_-$, the central term obeys the cocycle condition 
$$\chi_{m,n}^{q}=0 \qquad {\text {for}} \qquad\vert m-n\vert >2g_0.
\tag 26$$
\par
Finally in the $q\to 1$ limit the above algebra reduces to the
usual {\cd KN} algebra. Indeed, from (7) and (15b) one can ascertain
that ${\Cal L}^{\alpha}_{m}\overset {q \to 1}\to \longrightarrow
{\Cal L}_m$ and the operator coefficients in (21) reduce to (4a) 
while the central term becomes
$$\chi_{m,n}^{q}\overset {q \to 1}\to \longrightarrow \chi_{m,n}
= {{1}\over{12}}\sum_{k=0}^{2g_0-m-n} e_{m,k}^
+e_{n,2g_0-m-n-k}^+(m-g_0+k-1)(m-g_0+k)(m-g_0+k+1)\tag 27$$
which is equivalent to (4b) if the latter is evaluated using the
basis vector fields (2a).
\eject
\centerline {\sec References}
\vskip .3cm
\ref {1} Drinfeld V G  {\it Proc. Intern. Congress of
Mathematicians} Vol. 1 (Berkley, CA:University of California Press),1986,
 p 798.
\vskip.3cm
\ref {2} Curtright T L and Zachos C K  {\it Phys.
Lett.} B {\bf 243} 237 (1990).
\vskip.3cm
\ref {3} Chaichian M, Ellinas D and Popowicz Z 
{\it Phys. Lett.} B {\bf 248} 95 (1990).
\vskip.3cm  
\ref {4} Chaichian M and Pre${\check{\text s}}$najder P  
{\it Phys. Lett.} B {\bf 277} 109 (1992). 
\vskip.3cm
\ref {5} Polychronakos A P  {\it Phys. Lett.} B {\bf 256} 35 (1991).
\vskip.3cm
\ref {6} Devchand C and Saveliev M V  {\it Phys. Lett.} B
{\bf 258} 364 (1991).
\vskip.3cm
\ref {7} Aizawa N and Sato H  {\it Phys. Lett.} B {\bf 256} 185 (1991).
\vskip.3cm
\ref {8} Oh C H and Singh K  {\it J. Phys. A: Math. Gen.} 
{\bf 25} L149 (1992);  hep-th/9407081.
\vskip.3cm
\ref {9} Oh C H and Singh K  {\it Int. J. Mod. Phys.} {\bf 3A} 
(Proc. Suppl.) 383 (1993).
\vskip.3cm
\ref {10} Oh C H and Singh K  {\it J. Phys. A: Math. Gen.} 
{\bf 27} 3439 (1994);   hep-th/9408001.
\vskip.3cm
\ref {11} Belov A A and Chaltkian K D {\it Mod. Phys. Lett.} A 
{\bf 8} 1233 (1993).
\vskip.3cm
\ref {12} Krichever I M and Novikov S P  {\it Funct. Anal. Pril.}
{\bf 21} 46 (1987). 
\vskip.3cm
\ref {13} Krichever I M and Novikov S P  {\it Funct. Anal. Pril.}
{\bf 21} 47 (1987).
\vskip.3cm
\ref {14} Krichever I M and Novikov S P  {\it Funct. Anal. Pril.}
{\bf 23} 24 (1989).  
\vskip.3cm
\ref {15} Eguchi T and Ooguri H  {\it Nucl. Phys.} B {\bf 282} 308
(1987).   
\vskip.3cm
\ref {16} Mezincescu L, Nepomechie R I
{\it Nucl. Phys.} B {\bf 315} 43 (1989).
\vskip.3cm
\ref {17} Huang C S and Zhao Z Y {\it Phys. Lett.} B {\bf 220} 87 (1989).
\vskip.3cm
\ref {18} Chandrasekharan K  {\it Elliptic Functions},Springer
Verlag, 1984.
\vskip.3cm
\bye